\documentclass[10pt]{amsart}
\usepackage{amsmath, amssymb, latexsym}

\newtheorem{thm}[equation]{Theorem}
\newtheorem{pro}[equation]{Proposition}
\newtheorem{cor}[equation]{Corollary}
\newtheorem{lem}[equation]{Lemma}

\makeatletter\@addtoreset{equation}{section}\makeatother

\theoremstyle{definition}
\newtheorem{exa}[equation]{Example}

\newcommand{\sands}{\quad\text{and}\quad}
\newcommand{\spandsp}{\qquad\text{and}\qquad}
\newcommand{\co}{\colon}
\newcommand{\coo}{\,\colon\,}
\newcommand{\rta}{\rightarrow}
\renewcommand{\o}{\otimes}
\renewcommand{\=}{\,=\,}
\newcommand{\+}{\,+\,}
\newcommand{\isom}{\cong}
\newcommand{\pr}[2]{\ensuremath{\langle {#1,#2}\rangle}}
\newcommand{\prr}[2]{\ensuremath{\langle {#1\,,\,#2}\rangle}}
\newcommand{\sub}[1]{_{\text{${\scriptscriptstyle #1}$}}}
\newcommand{\subb}[2]{_{\text{${\scriptscriptstyle #1,#2}$}}}

\newcommand{\K}{\ensuremath{K}}   %the ring, or field

\newcommand{\s}{\sigma}
\newcommand{\cl}[1]{\ensuremath{c\ell (#1)}}
\newcommand{\clm}[1]{\ensuremath{c\ell_{\scriptscriptstyle M}(#1)}}
\newcommand{\cc}{\ensuremath{\mathcal C}}

\newcommand{\ci}{\ensuremath{\mathcal I}}
\newcommand{\cm}{\ensuremath{\mathcal M}}
\newcommand{\cu}{\ensuremath{\mathcal U}}
\newcommand{\frmod}[1]{\ensuremath{K\{{#1}\}}}      
\newcommand{\coal}[1]{\ensuremath{C(#1)}}           
\newcommand{\coaln}[2]{\ensuremath{C_{#1}(#2)}}     
\newcommand{\coalrk}[3]{\ensuremath{C_{#1,#2}(#3)}} 
\newcommand{\ic}[1]{\ensuremath{\widetilde{#1}}}    
\newcommand{\alg}[1]{\ensuremath{A(#1)}}            
\newcommand{\algn}[2]{\ensuremath{A_{#1}(#2)}}     
\newcommand{\algrk}[3]{\ensuremath{A_{#1,#2}(#3)}} 
\newcommand{\sect}[3]{\ensuremath{\binom{#1}{#2,#3}}}  
\newcommand{\msect}[3]{\ensuremath{\binom{#1}{#2,\dots,#3}}} 
\newcommand{\pow}[1]{\ensuremath{{\mathcal B}(#1)}}          
\newcommand{\powr}[2]{\ensuremath{{\mathcal B}_{#1}(#2)}}  
\newcommand{\wordn}[1]{\ensuremath{{\mathcal W}_{#1}}}     
\newcommand{\wordnr}[2]{\ensuremath{{\mathcal W}_{#1,#2}}} 
\newcommand{\word}{\ensuremath{\mathcal W}}  
\newcommand{\pos}[1]{\ensuremath{\pi (#1)}}     
\newcommand{\posk}[2]{\ensuremath{\pi_{#1}(#2)}}   
\newcommand{\chr}{\ensuremath{\text{\raisebox{.2ex}{$\chi$}}}}
\newcommand{\flag}{\ensuremath{(S_0,\dots,S_r)}}   
\newcommand{\sym}[1]{\ensuremath{\Sigma_{#1}}}    
\newcommand{\wmat}{\ensuremath{\lambda_{\scriptscriptstyle M}}}

\newcommand{\id}{\ensuremath{\iota}}         
\newcommand{\comp}{'}
\newcommand{\dual}{\ensuremath{^*}}
\newcommand{\pt}{\ensuremath{I}}           
\newcommand{\lp}{\ensuremath{Z}}           
\newcommand{\fr}[1]{\ensuremath{I_{#1}}}   
\newcommand{\ze}[1]{\ensuremath{Z_{#1}}}   
\newcommand{\mpt}[1]{\ensuremath{P_{#1}}}  
\newcommand{\cir}[1]{\ensuremath{C_{#1}}}  
\newcommand{\un}[2]{\ensuremath{U_{#1,#2}}}
\newcommand{\frc}{\ensuremath{\mathcal I}} 
\newcommand{\zec}{\ensuremath{\mathcal Z}} 
\newcommand{\all}{\ensuremath{\mathcal A}} 
\newcommand{\free}{\ensuremath{\mathcal F}}

\newcommand{\ws}[2]{\ensuremath{c(#1,#2)}}

\input{diagrams}
\begin{document}
\bibliographystyle{amsplain}

\title[A free subalgebra of the algebra of matroids]
{A free subalgebra of the algebra of matroids}
\author{Henry Crapo}
\author{William Schmitt}
\thanks{Schmitt partially supported by NSA grant 02G-134}
\email{crapo@ehess.fr and wschmitt@gwu.edu}

\begin{abstract}
This paper is an initial inquiry into the structure of
the Hopf algebra of matroids with restriction-contraction
coproduct.  Using a family of matroids introduced
by Crapo in 1965, we show that the subalgebra generated by a single point 
and a single loop in the dual of this Hopf algebra is free.
\end{abstract}

\maketitle

\section{introduction}

Major advances in combinatorial theory during recent decades rely upon
algebraic structures associated to combinatorial objects, and indeed,
often involve studies of combinatorial properties of algebraic systems
themselves.  In particular, Hopf algebras based on families of
combinatorial structures such as posets, graphs, permutations and
tableaux play an increasingly prominent role in contemporary
combinatorial theory and have been applied to a wide variety of
fields.  A major exception to this trend occurs in matroid theory,
where little  attention has been paid to naturally occurring
algebraic structures.  One such structure, introduced by one of the
present authors in \cite{sc:iha}, is a
Hopf algebra that may be associated to any family of matroids that is
closed under formation of minors and direct sums.  This Hopf algebra
has as basis the set of isomorphism classes of matroids belonging to
the given family, with product induced by the direct sum operation,
and coproduct of a matroid $M=M(S)$  given by $\sum\sub{A\subseteq
  S}M|A\o M/A$, where $M|A$ is the submatroid obtained by restriction
to $A$ and $M/A$ is the complementary contraction. A closely related
Hopf algebra was constructed by
Joni and Rota in
\cite{joro:cbc}, as the incidence coalgebra of a hereditary
family of geometric lattices.  In this case, attention
is restricted to simple matroids, and the subsets $A$ appearing
in the coproduct are taken to be flats.  These Hopf algebras were
also briefly considered in connection with the characteristic and
Tutte polynomials of matroids in \cite{korest:cft} and \cite{ku:mic}.

Similar constructions have arisen with increasing frequency in recent
years, as Hopf algebra techniques have been brought to bear on the
study of Feynman diagrams and renormalization processes in Physics
(\cite{cokr:rqf}, \cite{kr:oha}, \cite{brkr:rah}), Vassiliev's knot
invariants (\cite{chdula:vkiI}, \cite{chdula:vkiII},
\cite{chdula:vkiIII}, \cite{ko:vki}) and graph invariants
(\cite{elmo:nrm}, \cite{sa:hap}).  All of this  work has been
carried out in the context of graphs, which form an extremely
restricted class of matroids, and which have a grossly different
classification by isomorphism, save when attention is restricted to
3-connected graphs.
     
The present article is an initial inquiry into the structure of the
matroid Hopf algebra given in \cite{sc:iha}.  We prove that the
subalgebra of the dual algebra generated by ``point'' and ``loop''
(the two one-element matroids) is free.  (The question of whether of
not the corresponding subalgebra, in the context of graphs, is free,
which was posed by Lowell Abrams, remains open.)  We manage this proof
by restricting attention to a class of $2^n$ mutually nonisomorphic
matroids on an $n$ element set which we call ``freedom matroids''.
These matroids are obtained, starting from the empty matroid, by
successively adding points, at each stage either in a new dimension or
in general position in the top rank.  Freedom matroids were introduced
by the other present author, in \cite{cr:see}, in order to prove
that there are at least $2^n$ nonisomorphic matroids on $n$ elements.
The same matroids, presented as transversal matroids, were used in
\cite{we:bnm} to give a simplified proof of the same result.  Several
 characterizations of freedom matroids were given in
\cite{oxprro:mgs}, where it was also shown that the family of all
freedom matroids is closed under formation of minors and duals.  In
the present paper, we adduce a number of new combinatorial 
properties of freedom matroids.
This work is thus a useful adjunct to recent work
that has modeled these, and generalizations of these, matroids in
terms of Dyck paths (\cite{ar:tcm}) and lattice paths
(\cite{bodeno:lpm}), and other work, soon to appear (\cite{bode:lpm},
\cite{bogi:mpm}).

\section{Coalgebras of matroids}\label{sec:coalgebra}
Throughout this paper, we work over some commutative ring \K\ with
unit.  All modules, algebras and coalgebras are over \K, all maps
between such objects are assumed to be \K-linear, and all tensor
products are taken over \K.  Given any family of matroids \cm, we
write \ic\cm\ for the set of isomorphism classes of matroids belonging
to \cm, and denote by \frmod{\ic\cm} the free \K-module having \ic\cm\ 
as basis.  For any matroid $M=M(S)$, and $A\subseteq S$, we write
$M|A$ for the restriction of $M$ to $A$, and $M/A$ for the matroid on
$S\backslash A$ obtained by contracting $A$ from $M$.

The following result appeared in \cite{sc:iha}, as an example of the more
general construction of incidence Hopf algebras:

\begin{pro}\label{pro:coalgebra}
If \cm\ is a minor-closed family of matroids then \frmod{\ic\cm}
is a coalgebra, with coproduct $\delta$ and counit
$\epsilon$ determined by
$$
\delta (M)\= \sum_{A\subseteq S}M|A\o M/A
\spandsp
\epsilon (M)\= 
\begin{cases}
1, & \text{if $S=\emptyset$},\\
0, & \text{otherwise,}
\end{cases}
$$
for all $M=M(S)\in\cm$.  If, furthermore, the family \cm\ is 
closed under formation of direct sums, then \frmod{\ic\cm} is
a Hopf algebra, with product induced by direct sum.
\end{pro}
Whenever \cm\ is minor-closed, 
we shall write \coal\cm\ for the module \frmod{\ic\cm}
equipped with the above coalgebra structure.

We remark that in the statement of Proposition \ref{pro:coalgebra},
and in all that follows, we do not distinguish notationally between
matroids and their isomorphism classes; it will always be clear from
the context which is meant.
For the purposes of this article, we are 
interested primarily in the case in which
\cm\ is minor-closed and not necessarily closed under direct
sums and hence \coal\cm\ is only a coalgebra.
We don't give a complete proof of the proposition here, but
only mention that coassociativity of  $\delta$ follows 
directly from the basic identities 
$(M|T)|U=M|U$, $(M/U)/(T\backslash U)=M/T$ and $(M/U)|(T\backslash U)=
(M|T)/U$,
which hold for any matroid $M=M(S)$ and $U\subseteq T\subseteq S$.

In the case that \cm\ is closed under formation of direct
sums, a formula for the antipode of \cm\ may be deduced
from the formula for the antipode of an arbitrary incidence Hopf
algebra given in \cite{sc:iha}.

We will use the following notation for some specific matroids:
$$
  \begin{array}{ll}
\fr n = \un nn &\text{the {\it free matroid} of size $n$}\\
\ze n = \un 0n &\text{the {\it zero matroid} of size $n$}\\
\mpt n = \un 1 n &\text{the $n$-{\it point}}\\
\cir n = \un{n-1}n &\text{the $n$-{\it circuit}}\\
\pt = \fr 1 &\text{\it point}\\
\lp = \ze 1 &\text{\it loop},
\end{array}
$$
where, as usual, \un rn\ denotes the uniform matroid of 
rank $r$ on $n$ points.
\begin{exa}
  Let $L$ be the matroid shown in Figure \ref{fig:L},
consisting of points $a,b,c,d,e$ in
the plane, with $\{a,b,c\}$ and $\{a,d,e\}$ collinear.
If \cm\ is any minor-closed family containing $L$, then the
coproduct of $L$ in \coal\cm\ is given by
\begin{align*}
\delta (L) & \=
L\otimes\emptyset
\+ 4(\cir 3\oplus\pt)\o \lp\+\cir 4\o \lp
\+ 2\cir 3\o \mpt 2 
\+ 8\fr 3\o\ze 2\\
&\+ 6\fr 2\o (\mpt 2\oplus\lp)\+ 4\fr 2\o\mpt 3
\+ 4\pt\o N \+ \pt\o (\mpt 2\oplus\mpt 2)
\+ \emptyset\o L, 
\end{align*}
where $\oplus$ denotes the direct sum operation on matroids, and
$N$ is the three-point line with one of its points doubled.

\end{exa}

\begin{figure}[htbp]
\newcommand{\hllab}[1]
{\hspace{3.3ex}\text{\raisebox{0ex}{$\bf #1$}}}
\begin{diagram}[PostScript=Rokicki,abut,height=1.5em,width=1.5em,tight,thick]
&\hllab{c}&\bullet\\
&&\dLine\\
&\hllab{b}&\bullet\\
&&\dLine\\
&{\hspace{4.45ex}\text{\raisebox{-2.1ex}{$\bf a$}}}
&\bullet&\rLine&\bullet&\rLine
&\bullet&
{\hspace{-14.3ex}\text{\raisebox{-2.5ex}{$\bf d$\hspace{5.7ex}$\bf e$}}}
\end{diagram}
  \caption{}
  \label{fig:L}
\end{figure}

\begin{exa}\label{exa:free}
The family $\frc=\{\fr n\coo n\geq 0\}$ of all free matroids is
minor-closed, and the coalgebra \coal\frc\ is the free module $\K\{\fr
0,\fr 1,\dots\}$, with coproduct and counit given by $\delta (\fr
n)=\sum_{k=0}^n\binom nk\fr k\o\fr{n-k}$ and $\epsilon (\fr
n)=\delta_{n,0}$, for all $n\geq 0$.  Because \frc\ is also
closed under formation of direct sums, \coal\frc\ is in fact a Hopf
algebra.  Since \fr n\ is equal to the direct sum of $n$ copies
of \pt, we have $\fr n=\pt^n$ in \coal\frc, and thus \coal\frc\ is
the polynomial Hopf algebra $\K[\pt]$, with coproduct determined by
$\delta (\pt) = \pt\o 1 + 1\o \pt$.

Similarly, the family $\zec =\{\ze n\coo n\geq 0\}$ of all zero matroids
is closed under formation of minors and direct sums, and
\coal\zec\ is equal to the polynomial Hopf algebra $\K[\lp]$,
with $\delta (\lp) =\lp\o 1 + 1\o\lp$.
\end{exa}

Note that the coproducts in Example \ref{exa:free} are cocommutative.
This is because the operations of deletion and contraction on free and
zero matroids happen to coincide.  In fact, these are the only
matroids on which these operations coincide; if \cm\ is any
minor-closed family that contains matroids outside of $\frc\cup\zec$,
then the coalgebra \coal\cm\ is noncocommutative.

\begin{exa}\label{exa:uniform}
The class \cu\ of all uniform matroids
is minor-closed, and the coproduct on \coal\cu\ is
given by
$$
\delta (\un rn)\=\sum_{i=0}^r\binom ni\un ii\o\un{r-i}{n-i}
\+ \sum_{i=r+1}^n\binom ni\un ri\o\un 0{n-i},
$$
for all $n\geq r\geq 0$. 
If we adopt the convention
that $\un km =\un 0m$, for $k<0$ and $\un km=\un mm$, for
$k>m$, then the coproduct on \coal\cu\ takes the form
$$
\delta (\un rn)\=\sum_{i=0}^n\binom ni\un ri\o\un{r-i}{n-i},
$$
for all $n\geq r\geq 0$.
\end{exa}
\begin{exa}\label{exa:circuit}
  The subclass $\cc$ of \cu\ consisting of all circuits
and free matroids is minor-closed. The coalgebra \coal\cc\
is equal to $\K\{\fr 0,\fr 1,\dots,\cir 1,\cir 2,\dots\}$, with coproduct
determined by $\delta (\fr n)=\sum_{k=0}^n\binom nk\fr k\o\fr{n-k}$,
for $n\geq 0$, and $\delta (\cir m)=\cir m\o I_0 + \sum_{k=0}^{m-1}
\binom mk\fr k\o \cir{m-k}$, for all $m\geq 1$.
\end{exa}

Given a family \cm, and $n\geq 0$, we denote by $\cm_n$ the 
set of all matroids
belonging to \cm\ whose underlying sets have cardinality $n$;
and for $k,r\geq 0$, we denote by $\cm_{r,k}$ the set of
all matroids belonging to \cm\ that have rank $r$ and nullity $k$.
Writing \coaln n\cm\ and \coalrk rk\cm, respectively, for the free modules
\frmod{\ic\cm_n}\ and \frmod{\ic\cm_{r,k}}, we have
$$
\coal\cm\=\bigoplus_{n\geq 0}\coaln n\cm\=
\bigoplus_{r,k\geq 0}\coalrk rk\cm.
$$
\begin{pro}\label{pro:grade}
  If \cm\ is minor-closed, the families of submodules
 $\{\coaln n\cm\coo n\geq 0\}$ and $\{\coalrk rk\cm\coo r,k\geq 0\}$
 of \coal\cm, respectively, equip \coal\cm\ with the structure
of a graded, and bigraded, coalgebra.  If \cm\ is also closed
under formation of direct sums then \coal\cm\ is also thus
graded, and bigraded, as a Hopf algebra.
\end{pro}
\begin{proof}
    The first claim  follows immediately from the fact that, for
any matroid $M=M(S)$, and $A\subseteq S$, the rank of $M$
is equal to the sum of the ranks of $M|A$ and $M/A$, and
similarly for nullities. The second claim follows from the
fact that rank and nullity are additive functions with respect
to the disjoint sum operation on matroids.
  \end{proof}

  \begin{pro}\label{pro:dual}
If \cm\ is a minor-closed family and $\cm\dual =\{M\dual\coo M\in\cm\}$
then the map $D\sub\cm\co\coal\cm\rta\coal{\cm\dual}$, determined by
$M\mapsto M\dual$, for all $M\in\ic\cm$, is a coalgebra antiisomorphism.
In particular, if \cm\ is closed under duality, then $D\sub\cm$ is an
antiautomorphism of \coal\cm.
  \end{pro}
  \begin{proof}
    The map $D\sub\cm$ has inverse $D\sub{\cm\dual}$, and is thus
    bijective.  For any matroid $M=M(S)$, and $A\subseteq S$, we
    have the identities $(M|A)\dual = M\dual/(S\backslash A)$, and 
$(M/A)\dual=
    M\dual | (S\backslash A)$, from which it follows immediately that $\delta
    (D\sub\cm (M))=(D\sub\cm\o D\sub\cm)\cdot\tau\cdot \delta (M)$, where
    $\tau\co\coal\cm\o\coal\cm\rta\coal\cm\o\coal\cm$ is the twist
    map, determined by $M\o N\mapsto N\o M$, for all $M,N\in\cm$.
  \end{proof}

For all matroids $N_1$, $N_2$ and $M=M(S)$, the {\it section coefficient} 
\sect M{N_1}{N_2}\ is defined
as the number of subsets $A$ of $S$ such that $M|A\isom N_1$ and
$M/A\isom N_2$; hence if \cm\ is a minor-closed family, the coproduct
on \coal\cm\ is determined by
\begin{equation}
  \label{eq:seccop}
\delta (M)\=\sum_{N_1,N_2}\sect M{N_1}{N_2} N_1\o N_2,  
\end{equation}
 for all $M\in\cm$,
where the sum is taken over all (isomorphism classes of) matroids
 $N_1$ and $N_2$.
  We remark that there is no need
to restrict the sum in Equation \ref{eq:seccop} to matroids
$N_1$ and $N_2$ belonging to \cm; because the family \cm\ is
minor-closed, the section coefficient \sect M{N_1}{N_2}\ is
zero whenever $N_1$ or $N_2$ is outside of \cm.  
Another way of viewing this is the following:  If \all\ is
the class of all matroids, then the coproduct in \coal\all\
is given by Equation \ref{eq:seccop}; and if \cm\ 
is any minor-closed class then \coal\cm\ is a subcoalgebra
of \coal\all\ and thus the coproduct on \coal\cm\ is
given by the same formula as that for the coproduct on \coal\all.
\begin{exa}
Suppose that $M(S)$ is the matroid shown in Figure \ref{fig:bigmat}, and
that $N=\mpt 2\oplus\mpt 2$ is the matroid consisting of two double 
points. The section coefficient \sect M{\un 23}{N}\ is equal to one (rather
than two, as one might first guess) because, although 
there are two subsets $A$ of $S$ such that $M|A\isom\un 23$,
only for $A=\{a,b,c\}$ do we have $M/A\isom N$;
the contraction $M/\{a,d,e\}$, is a three point line with one
point doubled.
\end{exa}

\newcommand{\mlab}[1]
{\hspace{5.2ex}\text{\raisebox{1.ex}{$\bf #1$}}}
\newcommand{\mmlab}[1]
{\hspace{-6.2ex}\text{\raisebox{-0.9ex}{$\bf #1$}}}
 \newcommand{\nlab}[1]
{\hspace{6.4ex}\text{\raisebox{-2.3ex}{$\bf #1$}}}
\begin{figure}[ht]
\begin{diagram}[PostScript=Rokicki,abut,height=2.5em,width=2.3em,tight,thick]
%\begin{diagram}[abut,height=2.3em,width=2.3em,tight,thick]
&&\mlab{c}&\bullet&\rLine& &\bullet&\mmlab{g}\\
&\mlab{b}&\bullet\ldLine(1,1)&&&\ldLine(2,2)\\
\mlab{a}&\bullet\ldLine(1,1)&\rLine&&\bullet&\mmlab{f}\\
&\nlab{d}&\bullet\rdLine(1,1)&\\
&&\nlab{e}&\bullet\rdLine(1,1)&\\
\end{diagram}
\caption{}
\label{fig:bigmat}
\end{figure}

More generally, for matroids $N_1,\dots, N_k$ and 
$M=M(S)$, the {\it multisection coefficient} \msect M{N_1}{N_k} is
defined as the number of sequences $(S_0,\dots,S_k)$ such that
 $\emptyset =S_0\subseteq\cdots
\subseteq S_k=S$ and $(M|S_i)/S_{i-1}\isom N_i$, for $1\leq i\leq k$.
Hence the iterated coproduct $\delta^k\co\coal\cm\rta\coal\cm\o\cdots\o\coal\cm$
is determined by
$$
\delta^k (M)\=\sum_{N_1,\dots,N_k}\msect M{N_1}{N_k}
N_1\o\cdots\o N_k,
$$
for all $M\in\cm$.

\section{Algebras of matroids}\label{sec:algebras}

For any family of matroids \cm, we define a pairing $\pr\cdot\cdot\co
\frmod{\ic\cm}\times\frmod{\ic\cm}\rta\K$ by setting \pr MN\ equal to
the Kronecker delta $\delta\subb MN$, for all $M,N\in\cm$. This
pairing determines a pairing of $\frmod{\ic\cm}\o\frmod{\ic\cm}$ with
itself, by $\pr{M_1\o M_2}{N_1\o N_2} =
\pr{M_1}{N_1}\cdot\pr{M_2}{N_2}$, for all $M_1,M_2,N_1,N_2\in\cm$.  If
\cm\ is minor-closed, we may thus define a product on \frmod{\ic\cm},
dual to the coproduct on \coal\cm, by setting
\begin{equation}
  \label{eq:duality}
  \prr{N_1\cdot N_2}M\=\prr{N_1\o N_2}{\delta (M)}, 
\end{equation}
for all $M,N_1,N_2\in\cm$, thus making \frmod{\ic\cm}\ an
associative \K-algebra, with unit equal to the empty matroid.
We denote \frmod{\ic\cm}, equipped
with this algebra structure, by \alg\cm, and note that
\alg\cm\ is isomorphic to the graded dual algebra of
\coal\cm.

Writing \algn n\cm\ and \algrk rk\cm\ for the submodules
of \alg\cm\ generated, respectively, by matroids in \cm\
having $n$-elements, and those having rank $r$ and nullity $k$,
we have the direct sum decompositions:
$$
\alg\cm\=\bigoplus_{n\geq 0}\algn n\cm\=
\bigoplus_{r,k\geq 0}\algrk rk\cm,
$$
and it follows from Proposition \ref{pro:grade} that
\alg\cm\ is thus both a graded and bigraded algebra.
We also have the following result, dual to
Proposition \ref{pro:dual}.
\begin{pro}\label{pro:algdual}
If \cm\ is a minor-closed family and $\cm\dual =\{M\dual\coo M\in\cm\}$
then the map $D\co\alg\cm\rta\alg{\cm\dual}$, determined by
$M\mapsto M\dual$, for all $M\in\ic\cm$, is an algebra antiisomorphism.
In particular, if \cm\ is closed under duality, then $D$ is an
antiautomorphism of \alg\cm.
\end{pro}

By the definition of the pairing, the right-hand side of 
Equation \ref{eq:duality} is the coefficient of the basis 
element $N_1\o N_2$ in the coproduct $\delta (M)$ which, as noted in 
Equation \ref{eq:seccop}, is given by the section coefficient
\sect M{N_1}{N_2}. Since the left-hand side of \eqref{eq:duality}
is the coefficient of the basis element $M$ in the product
$N_1\cdot N_2$, it follows that
\begin{equation}
  \label{eq:secprod}
N_1\cdot N_2\=\sum_{M\in\ic\cm}\sect M{N_1}{N_2} M,
\end{equation}
for all $N_1,N_2\in\cm$. We emphasize that, in Equation
\ref{eq:secprod}, it is necessary to limit the summation to elements
of \ic\cm; because \coal\cm\ is a subcoalgebra of \coal\all, where
\all\ is the family of all matroids, it follows that \alg\cm\ is a
quotient of the algebra \alg\all.  Hence the product of $N_1$ and
$N_2$ in \alg\cm\ is the image of their product in \alg\all\ under the
projection homomorphism $\alg\all\rta\alg\cm$, which maps all matroids
$M\notin\ic\cm$ to zero.

\begin{exa}\label{exa:firstprods}
  Suppose that \cm\ is a minor-closed family containing
point \pt\ and loop \lp.
Then $\lp\cdot\pt = \pt\oplus\lp$ in \alg\cm.  If \cm\ contains
the double point \mpt 2 then
$\pt\cdot\lp=\pt\oplus\lp + 2\mpt 2$; otherwise,
$\pt\cdot\lp=\pt\oplus\lp$. If \cm\ contains the free matroid
\fr n\ then $\pt^n = n!\fr n$, and if \cm\ contains the zero 
matroid \ze n, we have $\lp^n = n!\ze n$ in \alg\cm.
\end{exa}
    
\begin{exa}
  Suppose that $L$ is the matroid shown in Figure \ref{fig:L} and
that $M$ is the matroid consisting of five points $a,b,c,d,e$ in
the plane, with $a,b,c$ collinear.   If \cm\ is any minor-closed
family that contains $L$, $M$ and the direct sum $\un 23\oplus\mpt2$
of the three-point line with a double point, then we have
$\un 23\cdot\mpt 2 = M + 2L + (\un23\oplus\mpt 2)$ in \alg\cm.
\end{exa}

\begin{exa}
If \cm\ contains the free matroid \fr r\ and zero matroid \ze k,
then the product $\fr r\cdot\ze k$ in \alg\cm\ is given by
$$
\fr r\cdot\ze k\= \sum (\text{$\#$ of bases of $M$})\cdot M,
$$
where the sum is over all matroids $M\in\ic\cm$ having rank $r$ and
nullity $k$.  On the other hand, for any $M\in\cm$ and $k\geq 0$, the
product $\ze k\cdot M$ is equal to $\binom{k+\ell}k\ze k\oplus M$,
where $\ell$ is the number of loops of 
$M$ if $\ze k\oplus M\in\cm$, 
and is equal to zero otherwise;
so in particular, $\ze
k\cdot\fr r = \ze k\oplus\fr r$ if \cm\ contains $\ze k\oplus\fr r$,
and $\ze k\cdot\fr r=0$, otherwise.
\end{exa}

\begin{exa}
Let \cc\ be the minor-closed family consisting of all
free matroids \fr n\ and circuits \cir k, for $n\geq 0$ 
and $k\geq 1$.  It follows from the coproduct formulas in Example 
\ref{exa:circuit} that the product in 
$\alg\cc\=\K[\fr 0,\fr 1,\dots,\cir 1, \cir2,\dots]$ is determined
by 
\begin{alignat*}{2}
\fr n\cdot\fr m& \=\binom{n+m}n\fr{n+m},&\qquad
\cir k\cdot\cir\ell & \=0,\\ \\
\fr n\cdot\cir k &\= \binom{n+k}n\cir{n+k},&\qquad 
\cir k\cdot\fr n&\=
\begin{cases}
  \cir k & {\text{if $n=0$,}}\\
  0 & {\text{otherwise,}}
\end{cases}
\end{alignat*}
for all $m,n\geq 0$ and $k,\ell\geq 1$.  
The dual family \cc\dual\ consists of all zero matroids \ze n\
and multiple points \mpt k, for $n\geq 0$ and $k\geq 1$.  By 
Proposition \ref{pro:algdual}, the product in \alg{\cc\dual}\
is determined by 
$\ze n\cdot\ze m=\binom{n+m}n\ze{n+m}$, 
$$
\mpt k\cdot\ze n \= \binom{n+k}n\mpt{n+k},
\quad \ze n\cdot\mpt k\=
\begin{cases}
  \mpt k & {\text{if $n=0$,}}\\
  0 & {\text{otherwise,}}
\end{cases}
$$
and $\mpt k\cdot\mpt\ell =0$, 
for all $m,n\geq 0$ and $k,\ell\geq 1$.  
\end{exa}

\section{Orderings of subsets and words}\label{sec:word}

For any set $S$ and $r\geq 0$, we denote by \pow S\ and \powr r S,
respectively, the set of all subsets and the set of all $r$-element
subsets of $S$.  In particular, for all $n\geq 0$, we write \pow n\ 
and \powr rn, respectively, for \pow{[n]}\ and \powr r{[n]}, where
$[n]$ denotes the set $\{1,\dots, n\}$.  Whenever we write a subset
of a linearly ordered set $S$ by listing its elements, we shall assume
that the list is written in the order induced by $S$; that is, if $S$
is linearly ordered, and $A=\{a_1,\dots,a_r\}\subseteq S$, then
$a_1<\cdots<a_r$ in $S$.  Throughout this paper we shall always assume
that $S$, whether linearly ordered or not, is a finite set.

For any linearly ordered $S$ and $r\geq 0$, we define a partial order
on \powr rS\ by setting $\{a_1,\dots, a_r\}\leq\{b_1,\dots,b_r\}$ if
and only if $a_i\leq b_i$ in $S$, for all $i\in [r]$.
Under this
ordering, \powr r S\ is a sublattice of the $r$-fold direct product of
linearly ordered sets $S\times\cdots\times S$, and is thus a
distributive lattice.  The Hasse diagram of \powr 2{\{a,b,c,d,e\}}\ is
shown in Figure \ref{fig:w52}.
\newcommand{\llab}[1]
{\hspace{-.5ex}\text{\raisebox{.2ex}{$\scriptstyle #1$}}}
\newcommand{\ali}[1] {\hspace{-9ex}\text{\raisebox{.2ex}{$\scriptstyle
      #1$}}} \newcommand{\alii}[1]
{\hspace{-11ex}\text{\raisebox{.2ex}{$\scriptstyle #1$}}}
\newcommand{\lllab}[1]
{\hspace{-3ex}\text{\raisebox{.2ex}{$\scriptstyle #1$}}}

\begin{figure}[htbp]
\mbox{\!\!\!\!\!\!\!\!
%\begin{diagram}[abut,height=1.2em,width=.8em,tight]
\begin{diagram}[PostScript=Rokicki,abut,height=1.2em,width=.8em,tight]
&&&&&&&&\llab{de}&\bullet\\
&&&&&&&&\ruLine(3,2)\\
&&&&&\llab{ce}& \bullet\\
&&&&&\ruLine(3,2)&&\luLine(3,2)\\
&&¦\llab{be}&\bullet&&&&&\llab{cd}&\bullet\\
&&\ruLine(3,2)&&\luLine(3,2)&&&&\ruLine(3,2)\\
\bullet&\ali{ae}&&&&\llab{bd}&\bullet\\
&\rdLine(3,2)&&&&\ldLine(3,2)&&\rdLine(3,2)\\
&&¦\llab{ad}&\bullet&&&&&\llab{bc}&\bullet\\
&&&&\rdLine(3,2)&&&&\ldLine(3,2)\\
&&&&&\llab{ac}&\bullet\\
&&&&&&&\rdLine(3,2)\\
&&&&&&&&\llab{ab}&\bullet\\
\end{diagram}}
\qquad\qquad\qquad\quad\mbox{
%\begin{diagram}[abut,height=1.2em,width=.8em,tight]
\begin{diagram}[PostScript=Rokicki,abut,height=1.2em,width=.8em,tight]
&&&&&&&&\lllab{00011}&\bullet\\
&&&&&&&&\ruLine(3,2)\\
&&&&&\lllab{00101}& \bullet\\
&&&&&\ruLine(3,2)&&\luLine(3,2)\\
&&¦\lllab{01001}&\bullet&&&&&\lllab{00110}&\bullet\\
&&\ruLine(3,2)&&\luLine(3,2)&&&&\ruLine(3,2)\\
\bullet&\alii{10001}&&&&\lllab{01010}&\bullet\\
&\rdLine(3,2)&&&&\ldLine(3,2)&&\rdLine(3,2)\\
&&¦\lllab{10010}&\bullet&&&&&\lllab{01100}&\bullet\\
&&&&\rdLine(3,2)&&&&\ldLine(3,2)\\
&&&&&\lllab{10100}&\bullet\\
&&&&&&&\rdLine(3,2)\\
&&&&&&&&\lllab{11000}&\bullet\\
\end{diagram}}\\
\caption{\text{The lattices $\powr 2{a,b,c,d,e}$ and
$\wordnr 5 2$}}\label{fig:w52}
\end{figure}

We extend the ordering on \powr rS\ to all
of \pow S\ by setting $B\geq A$ in \pow S\ if and only if
$B\geq A'$ in some \powr rS, for some subset $A'$ of $A$.
Hence, if $A=\{a_1,\dots,a_k\}$ and $B=\{b_1,\dots,b_r\}$, then
$A\leq B$ if and only if $r\leq k$ and $a_i\leq b_i$, for $1\leq i\leq r$.
Equipped with this ordering, \pow S\ is a distributive lattice that
contains each \powr rS\ as a sublattice.
\begin{lem}\label{lem:comp}
   For any linearly ordered set $S$, the map $\pow S\rta\pow S$ taking
$A\subseteq S$ to its complement in $S$ is a lattice
antiautomorphism.
\end{lem}
\begin{proof}
  Suppose that $A=\{a_1,\dots,a_k\}$ and $B=\{b_1,\dots,b_r\}$
are subsets of the linearly ordered set $S$ such that
$A\leq B$ in \pow S, that is, such that $r\leq k$
and $a_i\leq b_i$, for all $i\in [r]$.   If $A'=\{s_1,\dots,s_{n-k}\}$
and $B'=\{t_1,\dots,t_{n-r}\}$ are the complements of $A$ and $B$ in $S$,
then $n-r\geq n-k$, and
$s_j=j+|\{i\coo a_i<j\}|$ and $t_j=j+|\{i\coo b_i<j\}|$, for all $j$.
Since $a_i\leq b_i$, for all $i\in [r]$, it follows that
$|\{i\coo a_i<j\}|\geq |\{i\coo b_i<j\}|$, for all $j$.
Hence $s_j\geq t_j$, for $1\leq j\leq n-k$, and so
$A'\geq B'$ in \pow S.
\end{proof}

For any linearly ordered set $S$, we denote by $S_\varphi$ the {\it
  reversal\/} of $S$, that is, the set $S$ equipped with the opposite
ordering: $a\leq b$ in $S_\varphi$ if and only if $a\geq b$ in $S$.

\begin{lem}\label{lem:flip}
  For any linearly ordered set $S$, the identity map
is a lattice antiisomorphism $\powr rS\rta\powr r{S_\varphi}$.
\end{lem}
\begin{proof}
  It is immediate from the definition of the ordering
on \powr r S\ that $A\leq B$ in \powr rS\ if and only if
$A\geq B$ in \powr r{S_\varphi}.
\end{proof}

Given a word $w$ on the alphabet $\{0,1\}$, and $i\in\{0,1\}$, we
denote by $|w|_i$ the number of occurrences of the letter
$i$ in $w$.  For all $n\geq 0$, we write
\wordn n\ for the set of all words on $\{0,1\}$ having length $n$, and 
let $\wordnr n r= \{ w\in\wordn n \coo |w|_1=r\}$, for $0\leq r\leq n$.
For any linearly ordered set $S=\{e_1,\dots ,e_n\}$, let
$\chr\co\pow S\rta\wordn n$ be the function which maps $A\subseteq S$ to
the word $x_1\dots x_n$, where
$$
x_i \= \begin{cases}
1, & \text{if $e_i \in A$},\\
0, & \text{otherwise.}
\end{cases}
$$
Note that $\chr$ maps each \powr r S\ bijectively onto \wordnr n r\
and
that, under the natural identification of \wordn n\ with the set of
functions $S\rta\{0,1\}$, the function \chr\ simply
maps subsets of $S$ to their characteristic functions.

Define maps $\pi_k\co\wordnr nr\rta [n]$, for $1\leq k\leq r$, by
letting $\posk kw$ be the position of the $k${\small\it th} $1$ in 
$w\in\wordnr nr$.  
It follows that, for $S=\{e_1,\dots,e_n\}$, the map $\pi\co\
\wordnr nr\rta\powr r S$ which is inverse to $\chr$ is given by
$\pos w=\{e_{\posk 1w},\dots,e_{\posk rw}\}$, for all $w\in\wordnr nr$.
We define a
partial order on \wordnr nr\ by setting $v\leq w$ if and only if
$\posk kv\leq\posk kw$, for $1\leq k\leq r$.  
For example, the Hasse diagram of the lattice \wordnr 5 2\ is
given in Figure \ref{fig:w52}.
\begin{lem}\label{lem:latiso}
  For any linearly ordered set $S$, and $1\leq r\leq n=|S|$, the
map $\chr\co\powr r S\rta\wordnr nr$ is a lattice isomorphism.
\end{lem}
\begin{proof}
  It is immediate from the definition of $\chr$ that $A\leq B$
in \powr rS\ if and only if $\posk k{\chr (A)}\leq\posk k{\chr (B)}$,
for $1\leq k\leq r$.
\end{proof}
\begin{lem}\label{lem:wordineq}
  For all $v=x_1\cdots x_r$ and $w=y_1\cdots y_r$ in \wordnr nr,
the inequality $v\leq w$ holds if and only if $|x_1\cdots x_k|_1\geq
|y_1\cdots y_k|_1$, for $1\leq k\leq r$.
\end{lem}
\begin{proof}
The proof is immediate from the definitions.
\end{proof}

\section{Freedom Matroids}

By a {\it flag\/} on a finite set $S$ we shall mean a sequence \flag\ 
of subsets of $S$ such that $S_r=S$ and $S_{i-1}$ is a proper subset
of $S_i$, for $1\leq i\leq r$.  We do not require $S_0$ to
be empty.
\begin{pro}
For any flag \flag\ on a set $S$, the family 
$$
\ci\=\{ I\subseteq S\coo \text{$|I\cap S_i|\leq i$, for all $i$}\}
$$
is the collection of independent sets of a matroid $M\flag$,
of rank $r$, on $S$.
\end{pro}
\begin{proof}
  It is clear that \ci\ contains the empty set and is
closed under formation of subsets.  Now suppose that $I,J\in\ci$
with $|I|<|J|$.   If $|I\cap S_i|<i$ for all $i$, then for any
$x\in J\backslash I$ we have $|(I\cup x)\cap S_i|\leq i$ for all $i$, and
hence $I\cup x\in\ci$.  So we suppose that there exists some
$i$ such that $|I\cap S_i|=i$, and let $m$ be the maximal such
$i$.  Note that $m<r$, since $m=|I\cap S_m|\leq |I|<|J|=|J\cap S_r|\leq r$.

Now, since $|J\cap S_m|\leq m = |I\cap S_m|$, and $|J|>|I|$, we must have
$|J\cap S_m\comp|>|I\cap S_m\comp|$, where $S_m\comp$ denotes the complement of
$S_m$ in $S$, and hence the set $(J\backslash I)\cap S_m\comp$ is nonempty.
Let  $x$ be any element of $(J\backslash I)\cap S_m\comp$.
For $m<i\leq r$, we have $|I\cap S_i|<i$, and thus
$|(I\cup x)\cap S_i|\leq i$. Since $x\notin S_m$ we have
$(I\cup x)\cap S_i = I\cap S_i$, and so $|(I\cup x)\cap S_i|\leq i$, for
all $i\leq m$.  Thus $I\cup x\in\ci$.
\end{proof}

We refer to the matroid $M\flag$ as the {\it freedom matroid}
(see \cite{frfr})
defined by the flag \flag.  Note that it follows immediately
from the definition that each  $S_k$ is a flat of rank
$k$ in $M\flag$. 

If $M$ is a matroid on $S$ and $e\in S$, we denote by
$M\backslash e$ and $M/e$ the matroids obtained from $M$ by,
respectively, deleting and contracting $e$.  

\begin{pro}
For any freedom matroid $M=M\flag$ and $e\in S$, the 
deletion $M\backslash e$ and contraction $M/e$ are given by
$$
M\backslash e = M(T_0,\dots,T_r)\spandsp
M/e= M(T_0,\dots, T_{k-2},T_k,\dots,T_r),
$$
where $T_i=S_i\backslash e$, for all $i$, and $k=\min\{i\coo x\in S_i\}$.
\end{pro}
\begin{proof}
  The independent sets of $M\backslash e$ are the subsets of $S$ that
do not contain $e$ and contain no more than $i$ elements of each
$S_i$, which are precisely the independent subsets of
$M(T_0,\dots,T_r)$.

If $e$ is a loop in $M$, then $M/e=M\backslash e=M(T_0,\dots,T_r)$, which
agrees with the expression for $M/e$ given in the Proposition,
since $k=0$ in this case.  If $e$ is not a loop, then $A$ is independent
in $M/e$ if and only if $e\notin A$ and $A\cup e$ is independent in
$M$, that is $|(A\cup e)\cap S_i|\leq i$, for all $i$; in other words,
$|A\cap T_i|\leq i$, for $i<k$, and $|A\cap T_i|\leq
i-1$, for $i\geq k$.  Since $T_{k-1}\subseteq T_k$, the condition
$|A\cap T_k|\leq k-1$ implies that $|A\cap T_{k-1}|\leq k-1$ and hence
the latter inequality is redundant. Thus
$A$ is independent in $M/e$ if and only if $|A\cap T_i|\leq i$ for
$0\leq i\leq k-2$ and $|A\cap T_i|\leq i-1$, for $k\leq i\leq r$;
equivalently, if and only if $A$ is independent in
$M(T_0,\dots, T_{k-2},T_k,\dots,T_r)$.
\end{proof}
\begin{cor}[\cite{oxprro:mgs}]
  The class of freedom matroids is minor-closed.
\end{cor}

We now characterize the closure operators and closed
sets of freedom matroids.  We begin with the following proposition.

\begin{pro}\label{pro:closure}
The closure of an independent set $A$ in a freedom matroid
$M=M\flag$ is given by $\clm A = A\cup S_m$,
where $m=\max\{i\coo |A\cap S_i|=i\}$.
\end{pro}

\begin{proof}
  First note that $|A\cap S_0|=0$, because $A$ is independent, and thus
  such $m$ exists.  Now, since $|A\cap S_m|=m$, the set $A\cup x$ is
  dependent for all $x\in S_m\backslash A$, and thus $S_m\subseteq \clm A$.
On the other hand, for any $y\notin A\cup S_m$, the set
$A\cup y$ is independent, since
$|(A\cup y)\cap S_i|=|A\cap S_i|\leq i$, for $i\leq m$ and
$|(A\cup y)\cap S_i|\leq 1+|A\cap S_i|\leq i,$
for  $i>m$; hence $\clm A\subseteq A\cup S_m$.
\end{proof}
We may thus find the closure of an arbitrary set $A$
in a freedom matroid by applying Proposition \ref{pro:closure} to any 
maximal independent subset $B$ of $A$ and using the fact that
 $\cl B=\cl A$.
\begin{pro}\label{pro:closed}
A set $F\subseteq S$ is closed in $M\flag$ 
if and only if $F=A\cup S_m$,
for some $m\geq 0$ and $A\subseteq S\backslash S_m$ such that
$|A\cap S_i|<i-m$,
for all $i>m$; in which case the rank of $F$ is $m+|A|$.
\end{pro}
\begin{proof}
  Suppose that $F$ is closed and that $B$ is a basis for $F$.
By Proposition \ref{pro:closure}, $F=\cl B=B\cup S_m$ for some $m$
such that
$|B\cap S_m|=m$ and $|B\cap S_i|<i$, for all $i>m$.  Letting
$A=B\backslash S_m$, we thus have $F=A\cup S_m$ and $|A\cup S_i|<i-m$,
for all $i>m$.

On the other hand, suppose that $F=A\cup S_m$ for some $m\geq 0$
and $A\subseteq S\backslash S_m$, such that $|A\cup S_i|<i-m$, for all $i>m$.
Let $B$ be a basis for $S_m$.  Since $A$ is disjoint from
$S_m$, and thus also from $B$, and $|B|=m$, it follows from the
above inequality that $|(A\cup B)\cap S_i|\leq i$, for $i>m$,
and hence that $A\cup B$ is independent.  Since $m=\max\{i\coo
|(A\cup B)\cap S_i|=i\}$, it follows from Proposition 
\ref{pro:closure} that $A\cup S_m=\cl {A\cup B}$,
and is thus closed.
\end{proof}
Note that if we are given a closed set $F$ in $M\flag$, we can
express $F$ as $A\cup S_m$, according to Proposition 
\ref{pro:closed}, by letting $m=\max\{i\coo S_i\subseteq F\}$,
and taking $A=F\backslash S_m$.  
\begin{cor}\label{cor:flatsize}
If $F$ is any flat of rank $k$ in $M\flag$, then
$|F|\leq |S_k|$.  
\end{cor}
\begin{proof}
By Proposition \ref{pro:closed}, if $F$ is a flat of rank $k$ in $M\flag$
then $F=S_m\cup A$, for some $m$ and $A\subseteq S\backslash S_m$
with $|A|=k-m$.
Since $|S_k|-|S_m|\geq k-m$,
it follows that $|F|=|S_m|+|A|=|S_m|+k-m\leq |S_k|$.
\end{proof}

\section{Freedom matroids on ordered sets}

In the case that $S$ is linearly ordered it is convenient to consider
flags \flag\ such that each $S_i$ is an initial segment in the
ordering of $S$. In this case, the flag \flag\ is determined by
 $S$ together with the set
$\{1+\max{S_i}\coo 0\leq i\leq r-1\}$.  Hence if $S$ is linearly ordered
and we are given a subset $T=\{t_1,\dots,t_r\}$ of $S$, we may obtain
a flag $(T_0,\dots,T_r)$ on $S$ by setting $T_r=S$ and $T_i=\{s\in
S\coo s<t_{i+1}\}$, for $0\leq i\leq r-1$. 
We denote the freedom matroid
$M(T_0,\dots,T_r)$ by $M\sub T(S)$, or simply $M\sub T$, when the set
$S$ is understood.  If $T\subseteq [n]$ and $S=\{e_1,\dots,e_n\}$,
we also
write $M\sub T(S)$ for the matroid $M_{\alpha (T)}(S)$, where
$\alpha\co\pow n\rta\pow S$ is the natural bijection $i\mapsto e_i$.

\begin{pro}\label{pro:indep}
  If $S$ is linearly ordered and $T\subseteq S$, then
the family of independent sets of $M\sub T = M\sub T(S)$
is given by $\{ A\subseteq S \coo \text{$A\geq T$ in \pow S}\}$.
If $|T|=r$, then the family of bases of $M\sub T$ is given by
$\{ B\coo \text{$B\geq T$ in \powr rS}\}$.
\end{pro}
\begin{proof}
  Suppose that $T=\{t_1,\dots,t_r\}$
and $A=\{a_1,\dots,a_k\}$ in \pow S.  Since $T_r=S$, we have
$A=A\cap T_r$, and thus $|A\cap T_r|\leq r$ if and only if
$k\leq r$.  
Now for $0\leq i \leq r$, we have 
$A\cap T_i = \{ a_j\in A\coo \text{$a_j< t_{i+1}$ in $S$}\}$; therefore,
since $a_1<\cdots <a_k$ and $t_1<\cdots <t_r$, it follows
that $|A\cap T_i|\leq i$ if and only if $a_{i+1}\geq t_{i+1}$.
Hence $A$ is independent in $M\sub T$ if and only if $A\geq T$ in \pow S.
\end{proof}
\begin{exa}
  Suppose that $S=\{a,b,c,d,e,f,g\}$ and $T=\{b,e,f\}$. 
Then $M\sub T=M(T_0,T_1,T_2,T_3)$, where $T_0=\{a\}$, $T_1=\{
a,b,c,d\}$, $T_2=\{a,b,c,d,e\}$ and $T_3=S$.  The bases of
$M\sub T$ are the sets  $\{b,e,f\}$, $\{c,e,f\}$, $\{d,e,f\}$,
$\{b,e,g\}$, $\{c,e,g\}$, $\{d,e,g\}$, $\{b,f,g\}$,
$\{c,f,g\}$, $\{d,f,g\}$ and $\{e,f,g\}$.
\end{exa}
\begin{pro}
  For any linearly ordered $S$, and $T\subseteq S$,
the dual $M\sub T(S)\dual$ of the matroid $M\sub T(S)$ is equal
to $M\sub{T\comp}(S_\varphi)$, where
$T\comp$ is the complement of $T$ in $S$ and $S_\varphi$ 
is the reversal of $S$.  In particular, the
class of freedom matroids is closed under duality.
\end{pro}
\begin{proof}
Suppose that $|S|=n$ and $|T|=r$. It follows from
Proposition \ref{pro:indep} that the set of bases of
$M\sub T(S)\dual$ is given by $\{ B\comp\coo\text{$B\geq T$ in \powr rS}\}$,
which, according to Lemma \ref{lem:comp}, is equal to 
$\{ C\coo\text{$C\leq T\comp$ in \powr{n-r}S}\}$.
By Lemma \ref{lem:flip}, we have
$C\leq T\comp$ in \powr{n-r}S if and only 
$C\geq T\comp$ in \powr{n-r}{S_\varphi}, and hence
the result follows from Proposition \ref{pro:indep}.
\end{proof}

The following Lemma, which is a corollary of Proposition \ref{pro:indep},
will be used in the next section.
\begin{lem}\label{lem:up}
  Suppose that $M(S)=M\flag$ is a freedom matroid, where
$S$ is linearly ordered and each $S_i$ is an initial segment
in $S$, and let $A\subseteq S$ and $a\in A$.  If $b\in S\backslash A$ 
satisfies
$b>a$ in $S$, then $\rho ((A\backslash a)\cup b)\geq\rho (A)$.
\end{lem}
\begin{proof}
  Let $B$ be a maximal independent subset of $A$ that contains $a$.
Since $b>a$ in $S$, it follows that $(B\backslash a)\cup b>B$ in \pow S.
Hence, by Proposition \ref{pro:indep}, the set $(B\backslash a)\cup b$
is independent in $M$, and so $\rho ((A\backslash a)\cup b)\geq\rho (A)$.
\end{proof}

Recall from Section \ref{sec:word} that, given a word
$w\in\wordnr nr$, and $1\leq k\leq r$, we denote by $\posk kw$
the position of the $k${\small\it th} $1$ in $w$, and for
$S=\{e_1,\dots,e_n\}$, the bijection $\pi\co\wordnr nr\rta\powr r S$
is given by $\pos w=\{e_{\posk 1w},\dots,e_{\posk rw}\}$.  
We thus may define a mapping $w\mapsto M_w$ from \wordnr nr\
to the set of rank $r$ freedom matroids on $S$ by setting
$M_w=M_{\pos w}(S)$, for all $w\in\wordnr nr$.
\begin{exa}
  If $S=\{a,b,c,d,e,f,g,h,i,j,k,\ell\}$ and $w=001011001000$,
then $\pos w=\{c,e,f,i\}$. The sets
$S_i$ may be read off from the following table:
$$
\setcounter{MaxMatrixCols}{15}
\begin{matrix}
w: & 0&0&1&0&1&1&0&0&1&0&0&0\\
S_0: & a&b&\\
S_1: & a&b&c&d\\
S_2: & a&b&c&d&e\\
S_3: & a&b&c&d&e&f&g&h\\
S_4: & a&b&c&d&e&f&g&h&i&j&k&l,
\end{matrix}
$$
and $M_w=M_{\{c,e,f,i\}}$ is the freedom matroid
$M(S_0,S_1,S_2,S_3,S_4)$.
\end{exa}

When freedom matroids were first introduced, in \cite{cr:see}, they
were given the following recursive construction by single-element 
extensions: If $w$ is the empty
word, then $M_w$ is the empty matroid, and for $w=vx$, where $|x|=1$,
$M_w$ is obtained from $M_v$ as follows:\\
(i)  If $x=1$, add a point independently to $M_v$ in a new dimension,
that is, let $M_w=M_v\oplus I$;\\
(ii) If $x=0$, add a point $e$ to $M_v$ in general position in
the top rank, that is, let $M_w$ be the free extension of $M_v$ by $e$.
 
\begin{exa}
  If $w=001001010010$ and $S=\{a,b,c,d,e,f,g,h,i,j,k,l\}$, then
$M_w$ consists of loops $a$ and $b$, together with a triple point
$\{c,d,e\}$, collinear with distinct points $f$ and $g$, this line
being coplanar with general points $h$,$i$,$j$, with two additional
points $k$ and $l$ in general position in $3$-space.
\end{exa}

\section{Matroids and words}

Suppose that $M$ is a matroid of rank $r$ on an $n$-element set $S$,
having rank function $\rho$.  We associate to any maximal chain
$\emptyset=A_0\subset\cdots\subset A_n=S$ in the Boolean algebra $2^S$
the word $x_1\cdots x_n\in\wordnr nr$ defined by $x_i=\rho (A_i)-\rho
(A_{i-1})$, for all $i\in [n]$.  If the set $S=\{e_1,\dots,e_n\}$ is
linearly ordered, then there is a distinguished maximal chain
$A_0\subset\cdots \subset A_n$ in $2^S$, given by
$A_i=\{e_1,\dots,e_i\}$, for all $i\in [n]$.  The word $w\sub{M(S)}
=x_1\cdots x_n$ 
associated to this chain is thus determined by
$$
x_i \= 
\begin{cases}
0, & \text{if $e_i\in\cl{
\{e_1,\dots,e_{i-1}\}}$},\\
1, & \text{otherwise,}
\end{cases}
$$
for all $i\in [n]$. 
We refer to $w\sub{M(S)}$ as the {\it distinguished word} of $M(S)$.
Note that
$w\sub{M(S)}$ is also determined by the equality $|x_1\cdots x_i|_{1}=
\rho (\{e_1,\dots,e_i\})$, for all $i\in [n]$.
\begin{lem}\label{lem:dword}
  For any matroid $M(S)$ of rank $r$, with $S$ linearly ordered
of cardinality $n$,  the word $w=w\sub{M(S)}$ is determined by condition that
$\pi (w)=\min\{B\in\powr rS\coo\text{$B$ is 
a basis for $M$}\}$.
\end{lem}
\begin{proof}
Suppose $S=\{e_1,\dots,e_n\}$, and that the $1$'s in $w$ occur 
in positions $i_1,\dots,i_r$, 
so that $\pi (w)=\{e_{i_1},\dots,e_{i_r}\}$.  Since 
$e_{i_k}$ is not in the closure of $\{e_1,\dots,e_{i_k-1}\}$, 
for all $k\in [r]$, it follows that $\pi (w)$ is
independent, and thus is a basis for $M$.
If $B=\{b_1,\dots,b_r\}\subseteq S$ is such that $k\leq i_k$,
for some $k\in [r]$, then $\{b_1,\dots,b_k\}\subseteq
\{e_1,\dots,e_{i_k-1}\}$, which has rank $k-1$, and so $B$
is not a basis for $M$.  Hence any basis $B$ of $M$ satisfies
$B\geq\pi (w)$ in \powr rS.
\end{proof}

If $S=\{e_1,\dots,e_n\}$ is linearly ordered, then the symmetric group
\sym n\ acts naturally on $S$ by $\s (e_i)=e_{\s (i)}$, for all $i\in
[n]$, and thus we can identify \sym n\ with the group \sym S\ of
permutations of $S$.  For any $\s$ in \sym S\ (or in \sym n), we
denote by $S_\s$ the underlying set of $S$ equipped with the linear
order (or {\it reorder}) given by $\s (e_1)<\cdots <\s(e_n)$.  Hence,
$a\leq b$ in $S$ if and only if $\s (a)\leq\s (b)$ in $S_\s$, and so
$\s\co S\rta S_\s$ is a poset isomorphism.  The natural map $\pow
S\rta\pow{S_\s}$, given by $A\mapsto \s(A)$, for all $A\subseteq S$,
and also denoted by $\s$, is also a poset isomorphism.  We denote by
$\pi\sub\s$ the map $\wordnr nr\rta\powr r{S_\s}$, which takes a word
to the subset of $S_\s$ corresponding to positions of its $1$'s. Note
that $\pi\sub\s$ is equal to the composition $\s\pi$.

Given $A,B\subseteq S$ of equal cardinality, with complements $A'$ and
$B'$ in $[n]$, the {\it shuffle} $\s\sub{A,B}\in\sym S$ is the unique
permutation of $S$ which maps $B$ onto $A$, and thus also $B'$ onto
$A'$, whose restrictions to $B$ and $B'$ are order-preserving.  For
example, if $A=\{4,7\}$ and $B=\{1,5\}$ in $S=[7]$, then $\s\sub{A,B}=
4123756$ (where $\s =\s_1\cdots\s_n\in\sym n$ is the usual word
notation for permutations, indicating that $\s (i)=\s_i$, for all
$i$), or in cycle notation, $\s\sub{A,B}=(1432)(576)$.

\begin{lem}\label{lem:shuff}
   Suppose that $S$ is linearly ordered, and that
 $A\geq B$ in \pow S, where $|A|=|B|$, and 
let $\s=\s\sub{A,B}\in\sym S$ be the shuffle.
  If $C\subseteq S$ satisfies $C\geq A$ in \pow S, then $C\geq A$ in
  $\pow{S_\s}$.
\end{lem}

\begin{proof}
Suppose that the complements of $A=\{a_1,\cdots,a_r\}$ and 
$B=\{b_1,\dots,b_r\}$ in $S$ are $A'=\{a'_1,\cdots,a'_k\}$ and 
$B'=\{b'_1,\dots,b'_k\}$, respectively, so that the shuffle
$\s=\s\sub{A,B}$ is given by $b_i\mapsto a_i$ and $b'_j\mapsto
a'_j$, for all $i\in [r]$ and $j\in [k]$.  
Since $\s\co\pow S\rta\pow{S_\s}$ is an isomorphism,
it follows that
for any $C\subseteq S$, we have $C\geq A$ in \pow{S_\s}\ 
if and only if $\s^{-1}(C)\geq\s^{-1}(A)=B$ 
in \pow S.
Now suppose
that $C=\{c_1,\dots,c_m\}\geq A$ in \pow S, so that
$m\leq r$ and $c_i\geq a_i$, for all $i\in [m]$.
Since $A\geq B$ in
\pow S, it follows from Lemma \ref{lem:comp} that
$A'\leq B'$ in \pow S.  Hence $\s^{-1}(a)\leq a$, for all
$a\in A$, and $\s^{-1}(a')\geq a'$, for all $a'\in A'$.
Consider $c_i\in C$.  If $c_i\in A'$, then $\s^{-1}(c_i)\geq
c_i\geq a_i\geq b_i$. On the other hand, if $c_i\in A$, then
$c_i=a_j$, for some $j\geq i$ (since $c_i\geq a_i$), and so
$\s^{-1}(c_i)=\s^{-1}(a_j)=b_j\geq b_i$.  Hence $\s^{-1}(C)\geq
B$ in \pow S, and therefore $C\geq A$ in \pow{S_\s}.
\end{proof}

For any matroid $M(S)$ of rank $r$, where $S$ is linearly ordered of
cardinality $n$, we define a mapping $\wmat\co\sym S\rta\wordnr nr$
(or equivalently, $\wmat\co\sym n\rta\wordnr nr$) by setting $\wmat
(\s)=w\sub{M(S_\s)}$, for all
$\s\in\sym S$.  Note that, in particular, if $\id\in\sym S$ is the
identity permutation, then $\wmat (\id)=w\sub{M(S)}$ is the
distinguished word of $M(S)$.  We emphasize that the map $\wmat$
depends not only on the matroid $M=M(S)$, but on the linear ordering
of $S$.

For example, if $M$ is the matroid on $S=\{a,b,c,d,e,f,g\}$ shown
in Figure \ref{fig:bigmat}, and $\s\in\sym 7$ is the permutation
$6237154$, then $\wmat (\s)=1110010$.

\begin{pro}\label{pro:down}
Suppose that $M(S)$ is a rank $r$ matroid, with $S$ an $n$-element
linearly ordered set.
If $v\leq w\sub{M(S)}$ in \wordnr nr, then $\wmat (\s\sub{A,B})=v$,
where $A=\pi (w\sub{M(S)})$ and $B=\pi (v)$.
\end{pro}
\begin{proof}
By Lemma \ref{lem:dword}, $A=\pi (w\sub{M(S)})$ is the minimum
basis of $M$ in \powr rS.  Since $A\geq B=\pi (v)$ in
\powr rS, it follows from Lemma \ref{lem:shuff}
that $A$ is also the minimum basis of $M$ in \powr r{S_\s},
where $\s$ is the shuffle $\s\sub{A,B}$.  Since
$A=\s(B)=\s(\pi(v))=\pi\sub\s (v)$, it thus follows from
Lemma \ref{lem:dword} that $v=w\sub{M(S_\s)}$, that is, $\wmat (\s)=v$.
\end{proof}
\begin{cor}\label{cor:down}
  For any rank $r$ matroid $M$ on an $n$-element linearly ordered
set, the image of $\wmat$ is an order ideal in \wordnr nr.
\end{cor}
\begin{proof}
  The proof is immediate from Proposition \ref{pro:down}.
\end{proof}
It was shown in \cite{cr:see} (Theorem:
``Existence of a matroid with a given first word'') that in
the case in which $M=M_w$ is a freedom matroid, the word $w$ is
the maximum among words associated to $M$ by the map $\wmat$.
The following theorem is a strengthening of this result, giving a 
characterization of the words in the image of $\wmat$ whenever $M$
is a freedom matroid.
\begin{thm}  \label{thm:image}
If $M$ is the freedom matroid $M_w$ for some
$w\in\wordnr nr$, then the image of $\wmat\co\sym n\rta\wordnr nr$
is the principal order ideal $\{v\in\wordnr nr\coo v\leq w\}$.
\end{thm}
\begin{proof}
  Suppose that $M=M(S)=M_w$, where $S=\{e_1,\dots,e_n\}$ and
  $w=x_1\cdots x_n$ belongs to \wordnr nr.  It follows that
  $M=M\flag$, where $S_r=S$, and $S_{k-1}=\{e_1,\dots,e_{\posk
    kw-1}\}$, for $1\leq k\leq r$.  For any $\s\in\sym n$, the word
$\wmat (\s)=y_1\cdots y_n$ is determined by the condition
that $|y_1\cdots y_i|_1=\rho(\{e_{\s(1)},\dots,e_{\s(i)}\})$,
for $1\leq i\leq n$, and by Corollary \ref{cor:flatsize}, if
$\rho(\{e_{\s(1)},\dots,e_{\s(i)}\})=k$, for some $i$, then
$i\leq |S_k|=\posk{k+1}w-1$.  Since $\posk{k+1}w$ is the position
of the $(k+1)${\small\it st} one in $w$, it follows that 
$|x_1\cdots x_i|_1\leq k=|y_1\cdots y_i|_1$.  Hence, by
Lemma \ref{lem:wordineq}, we have $\wmat (\s)\leq w$.
The result thus follows from Corollary \ref{cor:down}.
\end{proof}

\begin{exa}
  Suppose that $M(S)=U_{2,4}\oplus P_2$ is the matroid consisting of
a four-point line and a double point. The image
of $\wmat$ in \wordnr 63\ (given any linear ordering on $S$) is
the order ideal $\{111000, 110100, 101100, 110010\}$, which has maximal
elements $110010$ and $101100$, and thus is not principal.
Hence, it follows from 
Theorem \ref{thm:image} that $M$ is not a freedom matroid.
\end{exa}

\begin{cor}[\cite{cr:see}]
There are precisely $2^n$ nonisomorphic freedom matroids
(and thus at least $2^n$ nonisomorphic matroids) on an 
$n$-element set.
\end{cor}
\begin{proof}
  Given a matroid $M$ on $S$, the definition of \wmat\ depends
on a choice of ordering of $S$, but the image of \wmat\ depends
only on the isomorphism class of $M$.  Hence, by Theorem \ref{thm:image},
if $v\neq w$, then the freedom matroids $M_v$ and $M_w$ are not
isomorphic.
\end{proof}

Recall that the {\it Bruhat order} (or {\it strong Bruhat order}) on \sym
n\ is determined by the condition that $\s$
covers $\tau =\tau_1\dots\tau_n$ in $\sym
n$ if and only if $\s$ may be obtained from $\tau $ by reversing a
single pair $(\tau_i,\tau_j)$, such that $i<j$ and $\tau_i<\tau_j$ 
and the number of inversions of $\s$ is one greater than the number
of inversions of $\tau$.  Under the assumptions $i<j$ and
$\tau_i<\tau_j$, the exchange $(\tau_i,\tau_j)$ increases the number
of inversions by one if and only if, for all $k$ with $i<k<j$, either
$\tau_k<\tau_i$ or $\tau_k>\tau_j$, which, in particular, is the case
if either $j=i+1$ or $\tau_j=\tau_i+1$.
For example, in the Bruhat order on \sym 4, 
the permutation $1423$ is covered by $4123$, $2413$ and $1432$.
Reversing the pair $(1,3)$ in $1423$ creates three new inversions, so that,
even though $3421$ is greater than $1423$, it is not a cover.
 The identity permutation is
the minimum element of \sym n, and the flip map $\varphi=n(n-1)\cdots 1$
is the maximum element.
\begin{pro}
  If $M=M_w$ for any $w\in\wordnr nr$, and $\sym n$ is  given
the Bruhat order, then $\wmat\co\sym n\rta\wordnr nr$ is an order-reversing
map.
\end{pro}
\begin{proof}
  Suppose that $M_w=M(S)= M\flag$, where
$S$ is linearly ordered and each $S_i$ is an initial segment in $S$.
Suppose that $\tau$ covers $\sigma$ in the Bruhat order on $\sym n$
and let $S_\s=\{e_1,\dots,e_n\}$ and $S_\tau=\{f_1,\dots,f_n\}$,
so that $e_k=f_k$ for all but two indices $i$ and $j$, where  
$$
i<j,\qquad e_i<e_j,\qquad f_j=e_i,\sands f_i=e_j.
$$
Letting $E_k=\{e_1,\dots,e_k\}$ and $F_k=\{f_1,\dots,f_k\}$, 
for all $k\in [n]$, we have $E_k=F_k$, for $1\leq k <i$ and $j<k\leq n$, and
since $e_j>e_i$ in $S$, it follows from Lemma \ref{lem:up} that 
$\rho (F_k)\geq\rho (E_k)$, for $i\leq k\leq j$.  
Letting $\wmat (\s)=x_1\cdots x_n$
and $\wmat (\tau )=y_1\cdots y_n$, we thus have
$|x_1\cdots x_k|_1=\rho (E_k)\leq\rho (F_k)=|y_1\cdots y_k|_1$, for
all $k\in [n]$, and hence $\wmat (\s)\geq \wmat (\tau)$,
by Lemma \ref{lem:wordineq}.
\end{proof}

\begin{exa}
  Suppose that $S=\{a,b,c,d\}$ and 
$M(S)=M\sub{0101}$, so that $a$ is a loop, $\{b,c\}$ a double point and
$d$ an isthmus in $M$.  The image of $\wmat\co\sym 4\rta\wordnr 42$ 
is the order ideal $\{1100, 0110, 1001, 1010\}$, and under 
$\wmat$, the
two permutations in the interval $[1234,1324]$ of \sym 4\ map
to $0101$, the four permutations in the interval $[1243,1432]$ map to 
$0110$, the four permutations in the interval $[2134,3214]$ map
to $1001$, the set $\{\s\coo\text{$\s\geq 2143$ and either
$\s\leq 3241$ or $\s\leq 4132$}\}$ maps to $1010$, and
the interval $[2413,4321]$ maps to $1100$.
\end{exa}

\section{The algebra of freedom matroids}

We now consider the algebra \alg\free\ corresponding to the
minor-closed class \free\ of freedom matroids.  Throughout
this section we shall assume that the ring \K\ is a field of characteristic
zero. 
The set $\{M_w\coo w\in\word\}$, where
$\word$ is the set of all words on
$\{0,1\}$, is a \K-vector space basis for \alg\free,
and the product is given by
$$
M_u\cdot M_v\=\sum_{w\in\word}\sect wuv M_w,
$$
where \sect wuv\ denotes the section coefficient \sect{M_w}{M_u}{M_v}.
As is the
case for any matroid algebra, $\alg\free$ is bigraded
by rank and nullity, and so $\alg\free
=\bigoplus_{r,k\geq 0}\algrk rk\free$,
where $\algrk rk\free$ has basis
$\{M_w\coo w\in\wordnr {r+k}r\}$, and the section
coefficient \sect wuv\ is zero whenever 
$w\notin \wordnr{|u|+|v|}{|u|_1+|v|_1}$.

In the proof of our main theorem below, we make use of the {\it
  incidence algebra} of the lattice \wordnr nr.  In general, the
incidence algebra $I(P)$ of a locally finite poset $P$ is
the \K-vector space of all functions $f\co P\times P\rta\K$ such
that $f(x,y)=0$, whenever $x\not\leq y$, equipped with the
{\it convolution} product:
$$
(fg)(x,z)\=\sum_{x\leq y\leq z}f(x,y)g(y,z),
$$
for all $f,g\in I(P)$, and $x\leq z$ in $P$.  
The convolution identity $\delta\in I(P)$
is given by $\delta (x,y)=\delta\sub{x,y}$, for all $x\leq y$ in
$P$.  An element $f\in I(P)$ is invertible if and only if $f(x,x)$
is a unit in \K, for all $x\in P$, in which case the convolution
inverse $f^{-1}$ is determined recursively by $f^{-1}(x,x)=
f(x,x)^{-1}$, for all $x\in P$, and 
$$
f^{-1}(x,z)\=f(z,z)^{-1}\sum_{x\leq y<z}f^{-1}(x,y)f(y,z)
\=f(x,x)^{-1}\sum_{x<y\leq z}f(x,y)f^{-1}(y,z),
$$
for all $x<z$ in $P$.

Recall that the matroids consisting of a single point and a single
loop are denoted by \pt\ and \lp, respectively, and note that
$\pt=M_1$ and $\lp =M_0$ are the freedom matroids corresponding to
words of length one.
\begin{thm}\label{thm:free}
The algebra \alg\free\ is free, generated by \pt\ and \lp.
\end{thm}
\begin{proof}
For any word $w=x_1\cdots x_n$ in \word, we denote by $P_w$ the
product $M_{x_1}\cdots M_{x_n}$ in \alg\free.
Since \alg\free\ is graded it suffices to show that
the set $\{P_w\coo w\in\wordnr nr\}$  is a basis for 
\algrk r{n-r}\free, for all $n\geq r\geq 0$. 
Given words $w,v\in\wordnr nr$, with $w=x_1\cdots x_n$,
we write \ws wv\ for the multisection coefficient
\msect v{x_1}{x_n}.  Observe that
\ws wv\ is equal to the number of permutations
$\s\in\sym n$ such that $\lambda\sub{M_v}(\s)=w$, and hence
Theorem \ref{thm:image} implies that \ws wv\ is nonzero if and
only if $w\leq v$ in the lattice ordering of \wordnr nr.  We
thus have
\begin{equation}\label{eq:mtop}
P_w\=\sum_{v\geq w}\ws wv M_v,
\end{equation}
for all $w\in\wordnr nr$,
where all coefficients are nonzero. Because $\ws wv=0$,
whenever $w\not\leq v$,  the function
$c$ belongs to the incidence algebra 
of \wordnr nr.  Since $\ws ww\neq 0$ for all $w$, and \K\ is a 
field of characteristic zero, it follows that $c$ has a convolution 
inverse $c^{-1}$, and therefore
$$
M_w\=\sum_{v\geq w}c^{-1}(w,v)P_v,
$$
for all $w\in\wordnr nr$. Hence the linear endomorphism of 
\algrk r{n-r}\free\ determined by $M_w\mapsto P_w$, for all 
$w\in\wordnr nr$, is invertible, and so $\{P_w\coo w\in\wordnr nr\}$ is
a basis for \algrk r{n-r}\free.
\end{proof}
Note that, since $P_v\cdot P_w=P_{vw}$ in \alg\free, for all
$v,w\in\word$, Theorem \ref{thm:free} can be restated as the fact
that the map $P_w\mapsto w$ defines an isomorphism from \alg\free\ 
onto the free algebra $\frmod\word=\K\langle\{0,1\}\rangle$, which has
concatenation of words as product.

The use of incidence algebras in the proof of Theorem
\ref{thm:free} can be avoided as follows:
Choose an ordering $w_1,\cdots,w_m$ of \wordnr nr\ such
that $i\leq j$, whenever $w_i\leq w_j$ in \wordnr nr\ 
(such as the opposite of lexicographic order)
and set $c_{ij}=c(w_i,w_j)$, for all $i\leq j$ in $[m]$.
Then $P_{w_i}\=\sum_{j=1}^m c_{ij}M_{w_j}$,
for all $i$,
 and by Theorem \ref{thm:image}, the matrix 
$C=(c_{ij})_{1\leq i,j\leq m}$
is upper-triangular, with nonzero entries along the main
diagonal.  Since $K$ is a characteristic zero field, 
$C$ is thus invertible, and hence the set $\{P_{w_i}\coo 1\leq i\leq m\}$
is a basis for \algrk r{n-r}\free.

\begin{cor}\label{cor:free}
  If \cm\ is 
any minor-closed family that contains the class \free\ of
freedom matroids, then the subalgebra of \alg\cm\ generated
by \pt\ and \lp\ is free.  
\end{cor}
\begin{proof}
For each word $w=x_1\cdots x_n\in\word$, let $Q_w$ denote
the product $M_{x_1}\cdots M_{x_n}$ in \alg\cm.
  Since $\free\subseteq\cm$, the algebra $\alg\free$ is
a quotient of $\alg\cm$, where the canonical homomorphism 
$\pi\co\alg\cm\rta\alg\free$ maps every freedom matroid in \cm\ 
to itself and every nonfreedom matroid to zero.
Since $\pi (Q_w)=P_w$, for all $w\in\word$ and, by Theorem
\ref{thm:free}, the $P_w$
are linearly independent in \alg\free, it follows that the
$Q_w$ are linearly independent in \alg\cm.  Hence
the subalgebra of \alg\cm\ generated by \pt\ and \lp\ is
free.
\end{proof}

\begin{exa}\label{exa:falg}
  If $S=\{a,b,c,d\}$, then the basis $\{M_w\coo w\in\wordnr 42\}$ 
of \algrk 22\free\ consists of the following
matroids:
$$
\begin{array}{ll}
M\sub{1100}\=\un 24  & \text{$a,b,c,d$ collinear}\\
M\sub{1010}  & 
\text{$\{a,b\}$ a double-point, collinear with points $c$ and $d$}\\
M\sub{1001}\=\mpt 3\oplus\pt  & 
\text{$\{a,b,c\}$ a triple-point, $d$ a distinct point}\\
M\sub{0110}\= \lp\oplus\un 23  & 
\text{$a$ a loop, $b,c,d$ collinear}\\
M\sub{0101}\=\pt\oplus\mpt 2\oplus\lp\;\;  & 
\text{$a$ a loop, $\{b,c\}$ a double-point, $d$ a distinct point}\\
M\sub{0011}\=\ze 2\oplus\fr 2  & \text{$a$ and $b$ loops, $c$ and $d$
distinct points}
\end{array}
$$
Listing $\wordnr 42$ in opposite lexicographic order, $\wordnr 42 =\{w_1,w_2,w_3,w_4,w_5,w_6\}=
\{1100, 1010, 1001, 0110, 0101, 0011\}$, the matrix $C$ of 
multisection coefficients $c_{ij}$ is given by
$$
\bordermatrix{
&{\raisebox{0ex}{$\scriptstyle 1100$}}&
{\raisebox{0ex}{$\scriptstyle 1010$}}&
{\raisebox{0ex}{$\scriptstyle 1001$}}&
{\raisebox{0ex}{$\scriptstyle 0110$}}&
{\raisebox{0ex}{$\scriptstyle 0101$}}&
{\raisebox{0ex}{$\scriptstyle 0011$}}\cr
{\scriptstyle 1100} & 24&20&12&12&8&4\cr
{\scriptstyle 1010} &  0 &4 &6 &6 &6&4\cr
{\scriptstyle 1001} &  0 &0 &6 &0 &4&4\cr
{\scriptstyle 0110} &  0 &0 &0 &6 &4&4\cr
{\scriptstyle 0101} &  0 &0 &0 &0 &2&4\cr
{\scriptstyle 0011} &  0 &0 &0 &0 &0&4\cr
}
$$\\
So, for example, $P_{1001}=\pt\cdot\lp\cdot\lp\cdot\pt$ is equal to
$6M_{1001}+4M_{0101}+4M_{0011}$ in $\alg\free$.  
Observe that $c_{34}$ is the only zero entry above the main diagonal $C$, 
which corresponds to the fact that $w_3=1001$ and $w_4=0110$
are the only two noncomparable elements of the lattice \wordnr 42.
Also note that, since the matrix entry $c(v,w)$ is equal to
the number of orderings of the underlying set of $M_w$
with corresponding  word equal to $v$, the
sum of the entries in each column of $C$ is equal to $4!$.
\end{exa}

\begin{exa}
Suppose that \cm\ is any minor-closed class containing
all freedom matroids and the smallest nonfreedom matroid
$D=P_2\oplus P_2$, consisting of two double-points, and let
$PL (\cm)$ be the subalgebra of \alg\cm\ generated by
\pt\ and \ze.  The matrix expressing the basis 
$\{Q_w\coo w\in\wordnr 42\}$ of 
$PL (\cm )\cap\algrk 22\cm$
in terms of the basis $\ic\cm_{2,2}=\{D\}\cup\{M_w\coo w\in\wordnr 42\}$ 
of \algrk 22\cm\ is given by
$$
\bordermatrix{
&{\raisebox{0ex}{$\scriptstyle 1100$}}&
{\raisebox{0ex}{$\scriptstyle 1010$}}&
 D&
{\raisebox{0ex}{$\scriptstyle 1001$}}&
{\raisebox{0ex}{$\scriptstyle 0110$}}&
{\raisebox{0ex}{$\scriptstyle 0101$}}&
{\raisebox{0ex}{$\scriptstyle 0011$}}\cr
{\scriptstyle 1100} & 24&20&16 &12&12&8&4\cr
{\scriptstyle 1010} &  0 &4& 8 &6 &6 &6&4\cr
{\scriptstyle 1001} &  0 &0& 0 &6 &0 &4&4\cr
{\scriptstyle 0110} &  0 &0& 0 &0 &6 &4&4\cr
{\scriptstyle 0101} &  0 &0& 0 &0 &0 &2&4\cr
{\scriptstyle 0011} &  0 &0& 0 &0 &0 &0&4\cr
}
$$\\
In this context, Corollary \ref{cor:free} amounts to
the observation that this matrix contains as a submatrix the nonsingular
matrix $C$ in the previous example, and thus
has independent rows.
\end{exa}

We now turn our attention to the coalgebra \coal\free\ of freedom matroids.
Recall from Section \ref{sec:coalgebra} that \coal\free\ has as
basis the set $\ic\free=\{M_w\coo w\in\word\}$ of all isomorphism classes
of freedom matroids,
and has coproduct determined by Equation \ref{eq:seccop}, so that
$$
\delta (M_w)\=\sum_{u,v\in\word}\sect wuv M_u\o M_v,
$$
for all $w\in\word$.  Hence if we define a coproduct
on the vector space \frmod\word, having all
$0$,$1$-words as basis, by $\delta (w)= \sum_{u,v}\sect wuv
u\o v$, then \frmod\word\ and \coal\free\ are isomorphic coalgebras
via the mapping $M_w\mapsto w$.
For example,
\begin{align*}
\delta (1010) & \= 1010\o \emptyset\+ 2\,(101\o 0)
\+ 2\,(110\o 0)\+ 10\o 10\\
&\+ 5\,(11\o 00)\+ 2\,(1\o 100)\+ 2\,(1\o 010)
\+ \emptyset\o 1010. 
\end{align*}
It is then an interesting exercise to give a description of
this coproduct solely in terms of the combinatorics of words.

Let $\{P'_w\coo w\in\word\}$ be the basis
of \coal\free\ which is dual to the basis $\{P_w\coo w\in\word\}$
of \alg\free\ via the pairing defined in the 
beginning of Section \ref{sec:algebras}, that is, such that
$\pr{P'_w}{P_v}=\delta\subb wv$,
for all $v,w\in\word$.
Equation \ref{eq:mtop} means that $\pr{M_v}{P_w}
=c(w,v)$, for all $v,w\in\word$, and so we have 
$$
  M_w \= \sum_{v\in\word}\pr{M_w}{P_v}P'_v
\= \sum_{v\leq w} c(v,w)P'_v
$$
for all $w\in\word$. Hence if $|w|=n$, and we write $\lambda$
for $\lambda\sub{M_w}$, we have
$$
M_w\=\sum_{\s\in\sym n}P'_{\lambda (\s)}.
$$
For example, referring to the matrix $C$ in Example
\ref{exa:falg}, we see that $M\sub{0110}=12P'\sub{1100} + 6P'\sub{1010}
+6P'\sub{0110}$ in \coal\free.

\begin{cor}\label{cor:cofree}
 The coalgebra \coal\free\ has basis $\{P'_w\coo w\in\word\}$
and coproduct given by
$$
\delta (P'_w)\=\sum_{uv=w}P'_u\o P'_v,
$$
for all $w\in\word$.
\end{cor}
\begin{proof}
  The result follows immediately from Theorem \ref{thm:free} by duality.
\end{proof}

Corollary \ref{cor:cofree} can be restated as saying that the map
determined by $P'_w\mapsto w$ is a coalgebra isomorphism from 
\coal\free\ onto the cofree coalgebra $\frmod\word$, which
has the deconcatenation coproduct $\delta (w)=\sum_{uv=w}u\o v$.

\end{document}